\documentclass[12pt,a4paper]{amsart}
\usepackage{amssymb,amscd}
\usepackage[margin=3cm]{geometry}

\newtheorem{thm}[subsection]{Theorem}
\newtheorem{lem}[subsection]{Lemma}
\newtheorem{prop}[subsection]{Proposition}
\newtheorem{cor}[subsection]{Corollary}
\theoremstyle{definition}

\newtheorem{example}[subsection]{Example}

\newtheorem{rmk}[subsection]{Remark}
\newtheorem{rmks}[subsection]{Remarks}

\numberwithin{equation}{subsection}

\newcommand{\N}{{\mathbb N}}
\newcommand{\Z}{{\mathbb Z}}
\newcommand{\Q}{{\mathbb Q}}

\newcommand{\End}{\operatorname{End}}

\newcommand{\Hom}{\operatorname{Hom}}

\newcommand{\ind}{\operatorname{ind}}

\newcommand{\GL}{\mathsf{GL}}

\newcommand{\T}{\mathsf{T}}
\newcommand{\M}{\mathsf{M}}
\newcommand{\D}{\mathsf{D}}
\newcommand{\gl}{\mathfrak{gl}}
\renewcommand{\sl}{\mathfrak{sl}}
\newcommand{\g}{\mathfrak{g}}

\newcommand{\divided}[2]{#1^{(#2)}}

\newcommand{\UU}{\mathbf{U}}

\newcommand{\A}{\mathcal{A}}
\newcommand{\bil}[2]{\langle #1, #2 \rangle}

\newcommand{\calF}{\mathcal{F}}

\newcommand{\vep}{\varepsilon}

\newcommand{\si}{\mathsf{i}}
\newcommand{\sj}{\mathsf{j}}
\newcommand{\E}{\mathsf{E}}
\newcommand{\e}{\mathsf{e}}
\newcommand{\row}{\mathrm{row}}
\newcommand{\col}{\mathrm{col}}
\newcommand{\U}{\mathfrak{U}}
\newcommand{\weight}{\operatorname{wt}}

\allowdisplaybreaks

\begin{document}
\title{A generic algebra associated to certain Hecke algebras} 

\author{Stephen Doty}
\address{Mathematics and Statistics, Loyola University Chicago, 
 Chicago, Illinois 60626 USA}
\email{doty@math.luc.edu}
\thanks{The first author was supported at Oxford by an EPSRC 
Visiting Fellowship;
all authors gratefully acknowledge support from Mathematisches 
Forschungsinstitut Oberwolfach, Research-in-Pairs Program.}

\author{Karin Erdmann}
\address{Mathematical Institute, 24-29 St Giles',
 Oxford OX1 3LB, UK}
\email{erdmann@maths.ox.ac.uk}

\author{Anne Henke}
\address{Mathematics and Computer Science,
 University of Leicester, University Road, 
 Leicester LE1 7RH, UK}
\email{A.Henke@mcs.le.ac.uk}

\begin{abstract} 
We initiate the systematic study of endomorphism algebras of
permutation modules and show they are obtainable by a descent from a
certain {\em generic} Hecke algebra, infinite-dimensional in general,
coming from the universal enveloping algebra of $\gl_n$ (or
$\sl_n$). The endomorphism algebras and the generic algebras are
cellular (in the latter case, of profinite type in the sense of
R.M. Green). We give several equivalent descriptions of these
algebras, find a number of explicit bases, and describe indexing sets
for their irreducible representations.
\end{abstract}
\maketitle

\section*{Introduction}
We study the intertwining spaces ${}_\lambda S(n,r)_\mu :=
\Hom_{\Sigma_r} (M^\mu, M^\lambda)$ between permutation modules
$M^\mu$, $M^\lambda$ for a symmetric group $\Sigma_r$. Here $\lambda$,
$\mu$ are given $n$-part compositions of $r$; that is, unordered
$n$-part partitions with $0$ allowed.  We are particularly interested
in the endomorphism algebras $S(\lambda) := {}_\lambda S(n,r)_\lambda$
of the modules $M^{\lambda}$.
These permutation modules are of central interest for the
representation theory of symmetric groups, and moreover provide a
natural link with the representation theory of general linear groups,
via Schur algebras.  Let $\E$ be a fixed $n$-dimensional vector space
over a field $K$. The symmetric group $\Sigma_r$ acts on the right on
$\E^{\otimes r}$ by place permutations. The endomorphism algebra
$$
S(n,r) = {\End}_{\Sigma_r}(\E^{\otimes r}) 
$$
is the Schur algebra. When $K$ is infinite its module category is
equivalent to the category of $r$-homogeneous polynomial
representations of $\GL_n(K)$.  As a module over $K\Sigma_r$, the
tensor space $\E^{\otimes r}$ is isomorphic to the direct sum $\oplus
M^{\lambda}$ where $M^{\lambda}$ is the transitive permutation module
corresponding to an $n$-part composition $\lambda$, that is the
$\Sigma_r$-orbit of weight $\lambda$ on the standard basis of
$\E^{\otimes r}$ (see \ref{ds}).  Hence $S=S(n,r)$ decomposes into a
direct sum of spaces $_\lambda S_\mu= {}_{\lambda}S(n,r)_{\mu}$, and
any such space is an $S(\lambda)$-$S(\mu)$ bimodule. We expect that a
systematic study of these algebras and bimodules will ultimately lead
to a better understanding of Schur algebras; however, the
finite-dimensional algebras $S(\lambda)$ are interesting in their own
right.  Moreover, the group algebra of $\Sigma_r$ appears in this
theory as the special case $S(\omega)={}_\omega S(r,r)_\omega$ for the
particular partition $\omega=(1^r)$, so our study may be regarded as
an extension of the study of symmetric groups to a broader context.

We also study, in the second part of the paper, certain
infinite-dimensional analogues of the $S(\lambda)$, ${}_\lambda
S_\mu$. These analogous objects occur naturally in Lusztig's
construction of the modified form $\dot{\U}$ of the universal
enveloping algebra $\U$ corresponding to the Lie algebra $\gl_n$ (over
base field $\Q$).  By definition $\dot{\U}$ is the direct sum of
spaces ${}_\lambda \U_\mu$ as $\lambda,\mu$ vary over the set $\Z^n$.
For given $\lambda \in \Z^n$, the space $\dot{\U}(\lambda):=
{}_\lambda \U_\lambda$ is an infinite-dimensional algebra, and if
$\lambda$ is an $n$-part composition then the finite-dimensional
algebra $S(\lambda)$ is a quotient of $\dot{\U}(\lambda)$. The space
${}_\lambda \U_\mu$ is a bimodule for
$\dot{\U}(\lambda)$-$\dot{\U}(\mu)$, and ${}_\lambda S_\mu$ is a
homomorphic image of ${}_\lambda \U_\mu$, by a linear map which is
compatible with the bimodule structures.  Similar statements apply if
one replaces $\dot{\U}(\gl_n)$ by $\dot{\U}(\sl_n)$.

Knowing the dimensions of the simple $S(\lambda)$-modules for every 
$n$-part composition $\lambda$ is equivalent to knowing the formal
characters of the simple $S(n,r)$-modules. Similarly,
knowing the dimensions of the simple $\dot{\U}(\lambda)$-modules for every 
$\lambda$ in $\Z^n$ is equivalent to knowing the formal
characters of the simple unital $\dot{\U}$-modules. 
See \S\ref{sec:Kostka}, \S\ref{sec:UKostka} for details.

\part{Endomorphism algebras of permutation modules}

We begin (in \S\S\ref{sec:wtsp}--\ref{sec:SA}) by summarizing basic
facts about weight spaces and Schur algebras, following ideas of J.A.\
Green and S.\ Donkin. Then (in \S\ref{sec:Hecke}) we show that the
intertwining space ${}_\lambda S(n,r)_\mu =
\Hom_{\Sigma_r}(M^\mu,M^\lambda) \simeq 1_\lambda S(n,r) 1_\mu$ is
isomorphic with the space
\begin{equation*}
\Hom_{\GL_n}(S^\mu \E, S^\lambda \E)
\end{equation*}
for $\lambda, \mu$ arbitrary $n$-part compositions of $r$, where
$S^\mu \E$ is a generalized symmetric power.
Here $1_{\lambda} \in S(n,r)$ corresponds to the projection onto
$M^{\lambda}$ with kernel $\oplus _{\mu \neq \lambda} M^{\lambda}$.
For $n\geq r$ this was proved in \cite{Donkin:SA2} but it holds in
general. In particular, it follows that $S(\lambda) \simeq
\End_{\GL_n}(S^\lambda \E)$.  We describe a connection between Kostka
duality and an indexing set for the simple $S(\lambda)$-modules (see
\S\ref{sec:Kostka}).

One can study Schur algebras (over $\Q$) as quotients of $\U$ and can
locate within $S_\Q(n,r)$ a certain integral form $S_\Z(n,r)$, a
quotient of a certain $\Z$-form $\U_\Z$ in $\U$ (see \ref{Ubas}). Then
$S_K(n,r) \simeq K\otimes _\Z S_\Z(n,r)$ for any field $K$. There are
similar integral forms for ${}_\lambda S_\mu$ and $S(\lambda)$, and we
give several natural $\Z$-bases for these spaces.

\section{Weight spaces} \label{sec:wtsp}

\subsection{} \label{wts}
Fix an arbitrary infinite field $K$.  We consider the affine algebraic
group $\GL_n = \GL_n(K)$ or the affine algebraic monoid $\M_n =
\M_n(K)$.  Let $\T_n \subset \GL_n$ (resp., $\D_n \subset \M_n$) denote
the subgroup (resp., submonoid) of diagonal matrices. Any rational left
$\GL_n$-module (or $\M_n$-module) $V$ has an eigenspace decomposition
relative to $\T_n$ (or $\D_n$): $V = \oplus_{\lambda\in \Z^n} \;
{}_\lambda V$ (resp., $V = \oplus_{\lambda\in \N^n} \; {}_\lambda V$),
where
\begin{equation}\label{wts:a}
{}_\lambda V = \{ v \in V \mid t \cdot v = t^\lambda v, \text{ all } t \in
\T_n\ (\text{or $\D_n$}) \}
\end{equation}
Similarly, if $V$ is a rational right $\GL_n$-module (resp., $\M_n$-module)
then $V$ has the decomposition: $V = \oplus_{\lambda\in \Z^n}\; V_\lambda$
(resp., $V = \oplus_{\lambda\in \N^n}\; V_\lambda$), where
\begin{equation}\label{wts:b}
V_\lambda = \{ v \in V \mid v \cdot t = t^\lambda v, \text{ all } t \in
\T_n \ (\text{resp., $\D_n$})  \}.
\end{equation}
Here $t^\lambda$ means $t_1^{\lambda_1} \cdots t_n^{\lambda_n}$, for $t =
\text{diag}(t_1, \dots, t_n)$.

Taking $V = K[\GL_n]$ (or $K[\M_n]$), regarded as bimodules for
$\GL_n$ (or $M_n$) via left and right translation of functions
($(g\cdot f)(x)=f(xg)$ and $(f\cdot g)(x)=f(gx)$, resp.), we
obtain direct sum decompositions of $K[\GL_n]$ (or $K[\M_n]$), both as
left and right modules. Let ${}_\lambda K$ (or $K_\lambda$) denote the
$1$-dimensional left (or right) $\T_n$- or $\D_n$-module afforded by
the character $\lambda$. The following (see \cite[3.5]{D:resolutions}
for details) is an immediate consequence of the definitions.

\begin{lem} \label{lem:ind}
  (a) For any $\lambda \in \Z^n$, the left (resp., right) weight space
  ${}_\lambda K[\GL_n]$ (resp., $K[\GL_n]_\lambda$) is isomorphic with the
  module $\ind_{\T_n}^{\GL_n} {}_\lambda K$ (or
  $\ind_{\T_n}^{\GL_n} K_\lambda$), as right (or left) rational
  $\GL_n$-modules.
  
 \noindent(b) For any $\lambda \in \N^n$, the left (or right)
  weight space ${}_\lambda K[\M_n]$ (resp., $K[\M_n]_\lambda$) is
  isomorphic with $\ind_{\D_n}^{\M_n} {}_\lambda K$ (or
  $\ind_{\D_n}^{\M_n} K_\lambda$) as right (or left) rational
  $\M_n$-modules.
\end{lem}

\subsection{} \label{act}
Although the decompositions in the preceding lemma are similar, it
should be noted that for $\lambda \in \N^n$ the modules
$K[\GL_n]_\lambda$, ${}_\lambda K[\GL_n]$ are infinite-dimensional
while $K[\M_n]_\lambda$, ${}_\lambda K[\M_n]$ are
finite-dimensional. Let $c_{ij}$ be the element of $K[\M_n]$ given by
evaluation of a matrix at its $(i,j)$th entry. Then $K[\M_n]$ may be
identified with the polynomial algebra $K[c_{ij}]$.  The algebra
$K[\GL_n]$ may be identified with the localization of $K[c_{ij}]$ at
the element $\det(c_{ij})$, and one may regard $K[\M_n]$ as a
subalgebra (in fact it is a sub-bialgebra) of $K[\GL_n]$ via the map
$K[\M_n] \to K[\GL_n]$ given by restricting functions from $\M_n$ to
$\GL_n$.  This gives a categorical isomorphism between polynomial
$\GL_n$-modules and rational $\M_n$-modules; this justifies our
interest in $\M_n$-modules. The bimodule structure on $K[\M_n]$ is
given explicitly by
\begin{equation}\label{act:a}\textstyle
g \cdot c_{ij} = \sum_k c_{ik} g_{kj}; \quad 
c_{ij} \cdot g = \sum_k g_{ik} c_{kj}  \qquad (g\in \M_n).
\end{equation}
These formulas also give the bimodule structure on $K[\GL_n]$, simply
by restricting $g$ to $\GL_n$. Moreover, by restricting to $\D_n$ we
obtain the equalities
\begin{equation}\label{act:b}
t \cdot c_{ij} =  c_{ij} t_{jj} = t^{\vep_j} c_{ij}; \quad 
c_{ij} \cdot t = t_{ii} c_{ij} = t^{\vep_i} c_{ij}  \qquad (t\in \D_n).
\end{equation}
This shows that $c_{ij}$ belongs to the weight space ${}_{\vep_j}
K[M_n]_{\vep_i}$; i.e., $c_{ij}$ has left (resp., right) weight
$\vep_j$ (resp., $\vep_i$), where $\vep_1, \dots, \vep_n$ is the
standard basis of $\Z^n$.

We have the following version of \cite[(2.7); (2.12)]{DW} or
\cite[proof of (3.4)(i)]{Donkin:MZ}, which provides an alternative
description of the weight spaces in $K[\M_n]$, in terms of symmetric
powers of the natural representation $\E$.

\begin{lem} \label{lem:sp}
  Let $\E$ be the $n$-dimensional vector space $K^n$, regarded as left
  (or right) $\M_n$-module by left (or right) matrix
  multiplication.  For any $\lambda \in \N^n$ there is an isomorphism
  between the induced module $\ind_{\D_n}^{\M_n} {}_\lambda K$ (resp.,
  $\ind_{\D_n}^{\M_n} K_\lambda$) and
  $$
  S^\lambda \E := (S^{\lambda_1} \E) \otimes \cdots \otimes
  (S^{\lambda_n} \E).
  $$
  as right (or left) rational $\M_n$-modules.
\end{lem}

\begin{proof} 
Let $(\e_i)$ be the canonical basis of $\E$. 
We may identify $S^a \E$ with homogeneous polynomials in the $\e_i$ of
degree $a$.  The left and right actions of $\M_n$ on $S^a \E$ are by
linear substitutions, and the natural action on $\E$.

By \eqref{act:b} the left weight space ${}_\lambda K[\M_n]$
is spanned by all monomials of the form
\begin{equation}\label{lem:sp:b}
\textstyle \prod_{i,j} c_{ij}^{a_{ij}} \quad( \sum_i a_{ij} = \lambda_j,
\text{ for all } j) 
\end{equation}
with $\M_n$ acting on the right as linear substitutions by the second
formula in \eqref{act:a} above.  Similarly, the right weight space
$K[\M_n]_\lambda$ is spanned by all monomials of the form
\begin{equation}\label{lem:sp:c}
\textstyle \prod_{i,j} c_{ij}^{a_{ij}} \quad( \sum_j a_{ij} = \lambda_i,
\text{ for all } i) 
\end{equation}
with $\M_n$ acting on the left as linear substitutions by the first
formula in \eqref{act:a} above.

The map taking an element of form \eqref{lem:sp:b} onto the element
\begin{equation}\label{lem:sp:e}
\textstyle (\prod_{i} \e_{i}^{a_{i1}})\otimes (\prod_{i}
\e_{i}^{a_{i2}})\otimes \cdots \otimes (\prod_{i} \e_{i}^{a_{in}})
\end{equation}
defines the desired isomorphism of right $\M_n$-modules. Similarly,
the map taking an element of form \eqref{lem:sp:c} onto the element
\begin{equation}\label{lem:sp:f}
\textstyle (\prod_{j} \e_{j}^{a_{1j}})\otimes (\prod_{j}
\e_{j}^{a_{2j}})\otimes \cdots \otimes (\prod_{j} \e_{j}^{a_{nj}})
\end{equation}
defines the desired isomorphism of left $\M_n$-modules.
\end{proof}

\begin{rmk}\label{rmk:1.5}
The proof shows that ${}_\lambda K[\M_n]$ (or $K[\M_n]_\lambda$) has a
basis in one-one correspondence with the set of $n\times n$ matrices
over $\N$ with column (resp., row) sums equal to $\lambda$.
\end{rmk}

\section{Schur algebras} \label{sec:SA}

\subsection{} \label{bas}
The Schur algebra can alternatively be constructed as the linear dual
of the coalgebra $A_K(n,r)$, the $K$-linear span of the
monomials in the $c_{ij}$ in $K[\M_n]$ of total degree $r$ (see
\cite[(2.4b)]{Green:book}).  This provides a basis for $S_K(n,r)$, as
follows. Given a pair of multi-indices $\si,\sj$ in $I(n,r)$, define
\begin{equation}\label{bas:b}
c_{\si,\sj} = c_{i_1 j_1} c_{i_2 j_2} \cdots c_{i_r j_r}. 
\end{equation}
Here $I(n,r)$ is the set of $\si = (i_1, \dots, i_r)$ where each $i_k$
belongs to $\{1, \dots, n\}$.  The commutativity of the variables
$c_{ij}$ implies that we have to take into account the equality rule
\begin{equation}\label{bas:c}
c_{\si,\sj} = c_{\si',\sj'} \Leftrightarrow \si' = \si\pi, \sj'= \sj\pi, 
\text{ some } \pi \in \Sigma_r,
\end{equation}
with respect to the obvious right action of $\Sigma_r$ on $I(n,r)$. As
a $K$-space, $A_K(n,r)^*$ has basis $\{ \xi_{\si,\sj} \}$ dual to the
basis $\{ c_{\si,\sj} \}$.  As with the $c_{\si,\sj}$ there is a
similar equality rule for the $\xi_{\si,\sj}$. The distinct
$\xi_{\si,\sj}$ provide the desired basis for the Schur algebra
$S_K(n,r)$.

In particular, by \cite[3.2]{Green:book} the distinct elements of the
form $\xi_{\si,\si}$ for $\si \in I(n,r)$ provide a set of orthogonal
idempotents which add up to the identity in $S_K(n,r)$. Given $\si$ we
set $\weight(\si) = (\lambda_1, \dots, \lambda_n) \in \N^n$ where
$\lambda_j$ is defined to be the number of $i_k$ which equal
$j$. We shall write
\begin{equation}\label{bas:d}
1_\lambda := \xi_{\si,\si} \quad (\lambda = \weight(\si) 
\text{ as above}) .  
\end{equation}
The distinct idempotents in this family are parametrized by the set
\begin{equation}\label{bas:e}
\Lambda(n,r) := \{ \lambda \in \N^n 
\mid \textstyle \sum_i \lambda_i = r  \}
\end{equation}
of $n$-part compositions of $r$. For any $S_K(n,r)$-module $V$, one
clearly has decompositions $V = \oplus_\lambda 1_\lambda V$ and $V =
\oplus_\lambda V 1_\lambda$. Moreover, one has by
\cite[3.2]{Green:book} identifications
\begin{equation} \label{bas:f}
1_\lambda V = {}_\lambda V, \qquad V 1_\lambda = V_\lambda
\end{equation}
with the weight spaces defined in \eqref{wts:a},
\eqref{wts:b}, for any $\lambda \in \Lambda(n,r)$.

\subsection{} \label{cntr}
It follows from \eqref{act:a} that $A_K(n,r)$, regarded as
$\M_n$-$\M_n$ bimodule, is $r$-homogeneous for either
action. Moreover, $A_K(n,r)$ is an $S_K(n,r)$-$S_K(n,r)$ bimodule (see
\cite[2.8]{Green:book}).  In particular, this means that the right
(resp., left) weight space $A_K(n,r)_\lambda$ (resp., ${}_\lambda
A_K(n,r)$) is a left (resp., right) $S_K(n,r)$-module.  For an
$S_K(n,r)$-module $V$, we denote by $V^\circ$ the contravariant dual
of $V$ (see \cite{Green:book}).

\begin{lem}\label{lem:A_K}
Let $\lambda \in \Lambda(n,r)$.
  
 \noindent(a) ${}_\lambda A_K(n,r) = {}_\lambda K[\M_n]$ and 
 $A_K(n,r)_\lambda = K[\M_n]_\lambda$.
  
 \noindent(b) There are isomorphisms $({}_\lambda A_K(n,r))^\circ 
 \simeq 1_\lambda S_K(n,r)$ (as right $S_K(n,r)$-modules) and
 $(A_K(n,r)_\lambda)^\circ \simeq S_K(n,r)1_\lambda$ (as left
 $S_K(n,r)$-modules).
\end{lem}

\begin{proof}\newcommand{\su}{\mathsf{u}} 
  (a) The space ${}_\lambda K[\M_n]$ is spanned by $\{c_{\si,\su} \mid
  \si \in I(n,r) \}$, where $\su$ is a fixed element of $I(n,r)$ of
  weight $\lambda$, so there is an inclusion ${}_\lambda K[\M_n]
  \subset {}_\lambda A_K(n,r)$. The opposite inclusion is
  obvious. This proves the first equality in (a), and the second is
  similar.
  
  (b) It is enough to prove the first statement, since the other case
  is similar.  One can check that the contravariant form $(\ , \ ):
  S_K(n,r) \times A_K(n,r) \to K$ defined in \cite[2.8]{Green:book},
  upon restriction to $1_\lambda S_K(n,r) \times {}_\lambda A_K(n,r)$,
  remains nonsingular.  In fact, $\{c_{\si,\su} \mid \si \in I(n,r)
  \}$ gives a basis for ${}_\lambda A_K(n,r)$ and $\{\xi_{\su,\sj}
  \mid \sj \in I(n,r) \}$ gives a basis for $1_\lambda S_K(n,r)$, and
  these two bases are dual relative to the contravariant form. The
  result now follows from \cite[2.7e]{Green:book}.
\end{proof}

\section{Hecke algebras associated to permutation modules}
\label{sec:Hecke}

We give a number of descriptions of the Hecke algebra $S_K(\lambda) =
\End_{\Sigma_r}(M^\lambda)$ and the bimodule ${}_\lambda S_K(n,r)_\mu =
\Hom_{\Sigma_r}(M^\mu, M^\lambda)$, for $\lambda,\mu \in \Lambda(n,r)$.

\subsection{} \label{ds} 
As right $K\Sigma_r$-modules, one has a decomposition
\begin{equation}\label{ds:a}
\E^{\otimes r} = \oplus_{\lambda}  M^\lambda.
\end{equation}
Here $M^\lambda$ is the transitive permutation module with basis all
$\e_{\si}=\e_{i_1}\otimes\cdots\otimes \e_{i_r}$ such that
$\weight(\si)=\lambda$.  The idempotents $1_{\lambda}$ of $S(n,r)$ as
in \eqref{bas:d} are precisely the projections corresponding to the
direct sum decomposition \eqref{ds:a}, see \cite{Green:book}.

By standard arguments we have an isomorphism of
$S_K(\lambda)$-$S_K(\mu)$ bimodules
\begin{equation}\label{ds:e}
1_\lambda S_K(n,r) 1_\mu \simeq \Hom_{\Sigma_r}(M^\mu, M^\lambda)
\end{equation}
and 
in particular an algebra isomorphism
\begin{equation}\label{ds:f}
1_\lambda S_K(n,r) 1_\lambda \simeq \End_{\Sigma_r}(M^\lambda).
\end{equation}
So we may identify $S_K(\lambda)$ and ${}_\lambda(S_K)_\mu$ with
$1_\lambda S_K(n,r) 1_\lambda$ and $1_\lambda S_K(n,r) 1_\mu$. This
shows in particular that $S_K(\lambda)$ is a type of Iwahori-Hecke
algebra; see \cite[6.1, Remark]{Green:book}. It also shows that,
whenever $n \ge r$, we have an isomorphism
\begin{equation}\label{ds:g}
1_\omega S_K(n,r)1_\omega \simeq K\Sigma_r,
\end{equation}
where $\omega \in \Lambda(n,r)$ is the special weight
\begin{equation}\label{ds:h}
\omega = (1,\dots,1, 0,\dots,0)  \qquad(\text{$r$ 1's}),
\end{equation}
since $M^{\omega}$ is the regular representation.  (See 
\cite[(6.1d)]{Green:book} for an explicit isomorphism.)

\medskip

One needs to consider $S_K(\lambda)$ and ${}_\lambda(S_K)_\mu$ only in
case $\lambda,\mu$ are dominant, that is $\lambda_1 \ge \lambda_2 \ge
\cdots \ge \lambda_n$, since these label the orbits of the Weyl group
$\Sigma_n$ acting on $\Lambda(n,r)$.

\begin{lem}\label{lem:iso} 
Let $\lambda, \mu \in \Lambda(n,r)$.

\noindent(a) For any $w \in \Sigma_n$, there is an algebra isomorphism
$S_K(\lambda) \simeq S_K(w \lambda)$.

\noindent(b) For any $w,w' \in \Sigma_n$, there is an isomorphism
${}_\lambda(S_K)_\mu \simeq {}_{w \lambda}(S_K)_{w' \mu}$ of
$S_K(\lambda)$-$S_K(\mu)$ bimodules.
\end{lem}

\begin{proof}
It is enough to show that $M^\lambda \simeq M^{w\lambda}$ for any $w
\in \Sigma_r$. But this is clear from the definition of $M^\lambda$
(see \ref{ds}).
\end{proof}

From the isomorphism \eqref{ds:e} and the multiplication rule
\cite[(2.3c)]{Green:book} we immediately obtain the following result,
which provides a basis for ${}_\lambda(S_K)_\mu$.

\begin{prop}\label{prop:Gbasis}
Let $\lambda, \mu \in \Lambda(n,r)$. A basis for $1_\lambda S_K(n,r)
1_\mu$ is given by the set of all $\xi_{\si,\sj}$ such that
$\weight(\si)=\lambda$, $\weight(\sj)=\mu$. In particular, the set of
all $\xi_{\si,\sj}$ ($\si,\sj \in I(n,r)$) satisfying $\weight(\si)=
\lambda = \weight(\sj)$ is a basis for $1_\lambda S_K(n,r) 1_\lambda
\simeq S_K(\lambda)$.
\end{prop}

In \cite{Green:JPAA}, J.A.\ Green introduced the codeterminant basis
of $S_K(n,r)$. There is a similar basis for $S_K(\lambda)$. Let $\nu
\in \Lambda(n,r)$ be a composition. Given a word $\si \in I(n,r)$ let
$T^\nu_\si$ denote the $\nu$-tableau obtained by writing the
components of $\si$ in order in the Young diagram of shape $\nu$.  The
weight of a tableau $T^\nu_\si$ is simply $\weight(\si)$. The
multi-index $\ell(\nu)$ is the word consisting of $\nu_1$ $1$'s,
followed by $\nu_2$ $2$'s, and so forth. The following is immediate
from the main result of \cite{Green:JPAA}.

\begin{prop}\label{prop:codet}
For $\lambda,\mu \in \Lambda(n,r)$, there is a basis for $1_\lambda
S_K(n,r) 1_\mu$ consisting of the codeterminants of the form
$Y_{\si,\sj}^\nu = \xi_{\si,\ell(\nu)} \xi_{\ell(\nu),\sj}$, such that
$\nu \in \Lambda(n,r)$ is dominant and $T^\nu_\si$, $T^\nu_\sj$ are
semistandard tableaux of weight $\lambda$, $\mu$, respectively.
\end{prop}

Thus the dimension of $S_K(\lambda)$ is the number of pairs of
semistandard tableaux (of some dominant shape $\nu$) of weight
$\lambda$.

\subsection{}\label{av} 
We can view $M^{\lambda}$ as a left $K\Sigma_r$-module. (Any right
$K\Sigma_r$-module can be viewed as a left $K\Sigma_r^{op}$-module,
but $K\Sigma_r^{op}$ is naturally isomorphic to $K\Sigma_r$, via the
involution induced by $\pi\to \pi^{-1}$ for $\pi \in \Sigma_r$.)

Assume that $n\ge r$. The Schur functor $\calF$ (see
\cite[6.1]{Green:book}) is the functor from left $S_K(n,r)$-modules to
left $K\Sigma_r$-modules defined by $V \to 1_\omega V$; it takes
$S^{\lambda}\E$ to $M^{\lambda}$. Indeed, by \ref{lem:ind} and
\ref{lem:sp} we have isomorphisms $\calF S^{\lambda}\E \simeq
{}_\omega K[\M_n]_\lambda \simeq (\E^{\otimes r})_\lambda =
M^\lambda$. 

Since $S^\lambda \E$ has a good filtration, the natural map
\begin{equation} \label{av:d}
\Hom_{S_K(n,r)}(S^\mu \E, S^\lambda \E) \to
\Hom_{\Sigma_r}(M^\mu,M^\lambda)
\end{equation}
is injective, for any $\lambda, \mu \in \Lambda(n,r)$. (See
\cite[1.8]{E} for details.) Thus, by comparing dimensions we see that
the map \eqref{av:d} is an isomorphism. (See \cite[2.4]{Donkin:SA2}.)
In particular, this gives an algebra isomorphism
\begin{equation}\label{av:e}
\End_{S_K(n,r)}(S^\lambda \E) \simeq \End_{\Sigma_r}(M^\lambda)
\end{equation}
for any $\lambda \in \Lambda(n,r)$.  Thus $S_K(\lambda) \simeq
\End_{S_K(n,r)}(S^\lambda \E)$ (provided $n \ge r$).

We would like to remove the restriction $n \ge r$ in the above.  So
suppose that $n < r$ and let $\E' = K^r$. There is another functor
taking $S_K(r,r)$-modules to $S_K(n,r)$-modules given by $M \to eM$
where $e$ is a certain idempotent in $S_K(r,r)$ (see
\cite[6.5]{Green:book}).  The idempotent $e$ is defined as follows.
We regard $\Lambda(n,r)$ as a subset of $\Lambda(r,r)$ via the
embedding $(\lambda_1, \dots, \lambda_n) \to (\lambda_1, \dots,
\lambda_n, 0, \dots, 0)$. Let $\Gamma$ be the image of this map; then
$e = \sum_{\beta \in \Gamma} 1_\beta$.  The functor takes $\E'$ to
$\E$ and $S^\lambda \E'$ to $S^\lambda \E$, for any $\lambda \in
\Lambda(n,r)$.  By Lemma \ref{lem:KE}, the functor induces a
surjection
\begin{equation}\label{av:f}
\Hom_{S_K(r,r)}(S^\mu \E', S^\lambda \E') \to
\Hom_{S_K(n,r)}(S^\mu \E, S^\lambda \E)
\end{equation}
for any $\lambda, \mu \in \Lambda(n,r)$.  We claim that this is an
isomorphism.  To see this, observe that from \ref{lem:ind},
\ref{lem:sp} and Frobenius reciprocity, the dimension of the
left-hand-side (right-hand-side) of \eqref{av:f} is $\dim_K {}_\lambda
K[\M_r]_\mu$ ($\dim_K {}_\lambda K[\M_n]_\mu$).  Since $\lambda, \mu$
are both in $\Lambda(n,r)$, embedded in $\Lambda(r,r)$ as above, we
see that these dimensions coincide. In fact, these dimensions are
given by the number of $r\times r$ matrices ($n\times n$ matrices)
over $\N$ with row sums $\mu$ and column sums $\lambda$; see Remark
\ref{rmk:1.5}.  Hence \eqref{av:f} is an isomorphism, as claimed.

In particular, this gives an algebra isomorphism
\begin{equation}\label{av:g}
\End_{S_K(r,r)}(S^\lambda \E') \simeq \End_{S_K(n,r)}(S^\lambda \E)
\end{equation}
for any $\lambda \in \Lambda(n,r)$.  We have removed the restriction $n
\ge r$ in \eqref{av:d} and \eqref{av:e} above.  In summary, we have
proved the following general result, with no restriction on $n,r$. 

\begin{prop}\label{pro:alt}
Let $n$ and $r$ be arbitrary positive integers, $\E = K^n$. 

\noindent(a) For $\lambda\in \Lambda(n,r)$, the algebra $S_K(\lambda)$
is isomorphic with
$$
\End_{S_K(n,r)}(S^\lambda \E) \simeq \End_{\GL_n}(S^\lambda \E)
$$

\noindent(b) For $\lambda, \mu \in \Lambda(n,r)$, the
$S_K(\lambda)$-$S_K(\mu)$ bimodule ${}_\lambda(S_K)_\mu$ is isomorphic
with
$$
\Hom_{S_K(n,r)}(S^\mu \E, S^\lambda \E) \simeq \Hom_{\GL_n}(S^\mu \E,
S^\lambda \E).
$$ 
\end{prop}

\begin{rmk}
(a) The above result, along with \eqref{ds:g}, proves that
$\End_{\GL_n}(\E^{\otimes r}) \simeq S_K(\omega) \simeq K\Sigma_r$,
provided $n \ge r$. 

(b) Since $\oplus_\mu S^\mu \E \simeq \oplus_\mu K[\M_n]_\mu =
A_K(n,r)$ (sums taken over all $\mu \in \Lambda(n,r)$), we see that
there is an algebra isomorphism $S_K(n,r) \simeq
\End_{\GL_n}(A_K(n,r))$. (This can be proved by other means.)
\end{rmk}

It remains to formulate the lemma refered to in the proof of the
proposition. In fact, we have a more general statement than what is
needed above.

\begin{lem}\label{lem:KE}
Let $S$ be a finite-dimensional algebra, $e \in S$ a nonzero idempotent. 
Suppose that $N$ is an injective $S$-module, and $eN$ is an injective
$eSe$-module. Then for any finite-dimensional $S$-module $M$ the 
restriction map 
$$
\Hom_S(M,N) \to \Hom_{eSe} (eM,eN) 
$$
is surjective. 
\end{lem}

\begin{proof}
The argument is similar to \cite[1.7]{E}.  We proceed by induction on
the composition length of $M$.  The inductive step is clear since the
functors ${\rm Hom}_S(-,N)$ and ${\rm Hom}_{eSe}(-, eN)$ are exact. So
we only have to prove surjectivity when $M$ is simple. Then either
$eM=0$, or $eM$ is simple and if $eM=0$ the statement follows
trivially. So assume now that $eM$ is simple. Clearly we may also
assume that $N$ is indecomposable injective.

We have a surjection $Se\otimes_{eSe}eM \to SeM$ and $SeM=M$ in this
case.  So we get an inclusion
\begin{equation}\label{eqn:KE}
0 \to \Hom_S(M, N) \to \Hom_S(Se\otimes_{eSe}eM, N)
\simeq \Hom_{eSe}(eM, eN)
\end{equation}
where the isomorphism comes from adjointness.
The general properties of adjointness give that
the composition of the inclusion and the isomorphism is
precisely the map induced by the functor $X\to eX$.

If $N$ is not the injective hull of $M$ then all the Hom spaces in
\eqref{eqn:KE} are zero. Otherwise they are 1-dimensional and then the
map of interest is an isomorphism (and hence is surjective).
\end{proof}

\section{Simple $S_K(\lambda)$-modules} \label{sec:Kostka}

There is a connection between Kostka duality and the problem of
parametrizing the simple $S_K(\lambda)$-modules.

\subsection{}
Let $\lambda \in \Lambda(n,r)$ be given. There is a functor
$\calF_\lambda$ from left $S_K(n,r)$-modules to left
$S_K(\lambda$-modules defined by $V \to 1_\lambda V$ (see
\cite[6.2]{Green:book}). The functor $\calF$ considered earlier is
$\calF_\omega$. For general $\lambda$, $\calF_\lambda$ is an exact
covariant functor mapping simple modules to simple modules or $0$.
Let $L(\mu)$ be the simple $S_K(n,r)$-module corresponding to a
dominant weight $\mu \in \Lambda(n,r)$.  By \cite[(6.2g)]{Green:book},
the collection of all nonzero $\calF_\lambda L(\mu)$ forms a complete
set of simple $S_K(\lambda)$-modules.

\subsection{}
For a dominant $\mu \in \Lambda(n,r)$, let $J(\mu)$ denote the
injective envelope of $L(\mu)$ in the category of rational
$\M_n$-modules. This is the contravariant dual of the projective cover
of $L(\mu)$, in the category of rational $\M_n$-modules (or the
category of $S_K(n,r)$-modules).  Since algebraic monoid induction
takes injectives to injectives, the generalized symmetric power
$S^\lambda \E \simeq \ind_{\D_n}^{\M_n} K_\lambda$ is injective as
$S(n,r)$-module, for any $\lambda \in \Lambda(n,r)$. Write $(S^\lambda
\E: J(\mu) )$ for the multiplicity of $J(\mu)$ in a Krull-Schmidt
decomposition of $S^\lambda \E$.  By Frobenius reciprocity one has
an equality
\begin{equation} \label{Kd}
( S^\lambda \E : J(\mu) ) = \dim_K 1_\lambda L(\mu) .
\end{equation}
(For more details see \cite{DW} or \cite[(3.4)]{Donkin:MZ}.)  The
equality is known as Kostka duality; the nonnegative integer in the
equality is the Kostka number, denoted by $\mathbf{K}_{\mu\lambda}$.
Note that $\mathbf{K}_{\mu\lambda}$ may be equivalently defined to be
the multiplicity of a Young module $Y^\mu$ in a Krull-Schmidt
decomposition of $M^\lambda$; see \cite[(3.5), (3.6)]{Donkin:MZ}.

\begin{prop} \label{prop:Kd}
Let $\lambda \in \Lambda(n,r)$ be fixed.  The isomorphism classes of
simple $S(\lambda)$-modules are the ${}_\lambda L(\mu) = 1_\lambda
L(\mu)$ for which $\mathbf{K}_{\mu\lambda} \ne 0$.
\end{prop}

\begin{proof} 
Combine \cite[(6.2g)]{Green:book} with Kostka duality.
\end{proof}

\begin{rmk}
Donkin \cite[p.~55, Remark]{Donkin:MZ} gives a more precise necessary
and sufficient condition on $\lambda,\mu$ for $\mathbf{K}_{\mu\lambda}
\ne 0$, obtained by combining Steinberg's tensor product theorem with
Suprunenko's theorem \cite{Su}.
\end{rmk}

\section{PBW basis} \label{sec:PBW}

In this section, we work with the Lie algebra $\gl_n$ and its
eveloping algebra $\U$ over the rational field $\Q$. We also consider
the algebra of distributions (the hyperalgebra) $\U_\Z$ of the
algebraic $\Z$-group scheme $\mathbf{GL}_n$ defined by $\Z[c_{ij};
(\det(c_{ij}))^{-1}]$.

\subsection{}\label{Ubas}
Let $e_{ij}$ be the $n\times n$ matrix whose unique nonzero entry is a
$1$ in the $(i,j)$th position.  Set $f_{ij} = e_{ji}$ and $H_i =
e_{ii}$.  The set
\begin{equation}\label{Ubas:a}
\{ f_{ij} \mid i<j\} \cup \{ H_i\} \cup \{ e_{ij} \mid i<j \}
\end{equation}
is a Chevalley basis for $\gl_n$. We regard these as elements of the
universal enveloping algebra $\U = \U(\gl_n)$.  
Note that $\U$ is generated by the elements
\begin{equation} \label{Ugen}
e_i = e_{i,i+1}, \  f_i = e_{i+1,i} \ \ (1 \le i \le n-1), \quad 
H_i\ \ (1\le i \le n).  
\end{equation}
For an element $X$ and an integer $a\ge 0$, set $\divided{X}{a} :=
X^a/(a!)$, $\binom{X}{a} := X(X-1)\cdots(X-a+1)/(a!)$.  The
hyperalgebra $\U_\Z$ of the $\Z$-group $\mathbf{GL}_n$ is the
$\Z$-subalgebra of $\U$ generated by all $\divided{f_i}{a}$,
$\binom{H_i}{b}$, $\divided{e_i}{c}$ ($a,b,c \ge 0$); see \cite[II,
Chapter 1]{Jantzen}. The set of all products of the form
\begin{equation}\label{Ubas:b}
\prod_{i<j} \divided{f_{ij}}{a_{ij}} \prod_i \binom{H_i}{b_i}
\prod_{i<j} \divided{e_{ij}}{c_{ij}}
\end{equation}
for nonnegative integers $a_{ij}$, $b_i$, $c_{ij}$, forms a $\Z$-basis
of $\U_\Z$, where the products among the $f$'s (resp., $e$'s) are
taken in some arbitrary, but fixed, order.  With similar conventions
on the order of products, the set
\begin{equation}\label{Ubas:c}
\prod_{i<j} \divided{e_{ij}}{c_{ij}} \prod_i \binom{H_i}{b_i} 
\prod_{i<j} \divided{f_{ij}}{a_{ij}}
\end{equation}
is another $\Z$-basis of $\U_\Z$.  

\subsection{} \label{bc}
By differentiating the homomorphism $\rho: \GL_n\to {\rm
End}_{\Sigma_r}(\E^{\otimes r})$ one obtains a representation $\gl_n
\to \End(\E^{\otimes r})$; this extends uniquely to a representation
\begin{equation}\label{Ubas:d}\begin{CD}
\U @>d\rho>> \End(\E^{\otimes r})
\end{CD}\end{equation}
and $S(n,r)$ (over $\Q$) is the image of this
representation. Moreover, $S_\Z(n,r)$ is the image of $\U_\Z$ under
the same map. In fact, one can consider the $\U_\Z$-invariant lattice
$\E_\Z$ (the $\Z$-span of the standard basis $(\e_i)$ on $\E$) and
check that the map $d\rho$ takes $\U_\Z$ into $\End(\E_\Z^{\otimes
r})$; then $S_\Z(n,r)$ is the image $d\rho(\U_\Z)$ \cite{CL}.  One has
isomorphisms
\begin{equation}\label{Ubas:e}
S_K(n,d) \simeq K \otimes_\Z S_\Z(n,r), \quad \U_K \simeq K \otimes_\Z
\U_\Z
\end{equation}
for a field $K$, where $\U_K$ is the hyperalgebra of the $K$-group
$\mathbf{GL}_n$. Thus $S_K(n,r)$ is a homomorphic image of $\U_K$.
Denote the images in $S(n,r)$ of the elements $\divided{f_{ij}}{a}$,
$\binom{H_i}{b}$, $\divided{e_{ij}}{c}$ by the same symbols.

\begin{lem} \label{lem:idem}
The element $1_\lambda$ in $S(n,r)$ coincides with the element given
by the product $\prod_i \binom{H_i}{\lambda_i}$, for any $\lambda \in
\Lambda(n,r)$.
\end{lem}

\begin{proof}
The element $1_\lambda$ was defined in \eqref{bas:d}.  It is known
(see \cite[3.2]{Green:book}) that, if $M$ is any $S(n,r)$-module, then
$1_\lambda M$ coincides with the $\lambda$-weight space of $M$ (for
$\lambda \in \Lambda(n,r)$), considered as left
$\GL_n$-module. Similarly, one can show that the element
$1_\lambda^\prime := \prod_i \binom{H_i}{\lambda_i}$ acts as zero on
all weight space except the $\lambda$ one, and acts as $1$ there. This
proves that $1_\lambda$ and $1_\lambda^\prime$ act the same on all
modules $M$. To finish, we need the existence of a faithful module,
since on such a module the difference $1_\lambda - 1_\lambda^\prime$
acts as zero, hence must equal zero in the algebra. Since $\E^{\otimes
r}$ is faithful, as $S(n,r)$-module, the proof is complete.
\end{proof}

\subsection{}\label{PSA}
\cite[Theorem 2.3]{DG:PSA} can be reformulated as follows.  Given an
$n\times n$ matrix $A=(a_{ij})$, let $\lambda^+(A) := (\lambda_1^+,
\dots, \lambda_n^+)$, where $\lambda_j^+ = a_{jj} + \sum_{i<j}
(a_{ij}+a_{ji})$ for each $j$.  Similarly, $\lambda^-(A) :=
(\lambda_1^-, \dots, \lambda_n^-)$, where $\lambda_j^- = a_{jj} +
\sum_{i>j} (a_{ij}+a_{ji})$ for each $j$.  Let $\Theta(n,r)$ be the
set of $n \times n$ matrices over $\N$ whose entries sum to $r$.  Then
by \cite[Theorem 2.3]{DG:PSA} the set of products of the form
\begin{equation} \label{PSA:a}
(\prod_{i<j} \divided{f_{ij}}{a_{ji}})\, 1_{\lambda^-(A)}\,
(\prod_{i<j} \divided{e_{ij}}{a_{ij}}) \quad (A \in \Theta(n,r))
\end{equation}
is a $\Z$-basis for $S_\Z(n,r)$.  Similarly, the set
\begin{equation}\label{PSA:b}
(\prod_{i<j} \divided{e_{ij}}{a_{ij}})\, 1_{\lambda^+(A)}\,
(\prod_{i<j} \divided{f_{ij}}{a_{ji}}) \quad (A \in \Theta(n,r))
\end{equation}
is another basis for $S_\Z(n,r)$. As usual, one has to take products
of $f$'s (resp., $e$'s) with respect to an arbitrary, but fixed,
order.  One can easily see that any element of the form \eqref{PSA:b}
has content (in the sense defined in \cite[p.~1911]{DG:PSA}) not
exceeding $\lambda^+(A)$; conversely, given a monomial of the form
$e_A 1_\mu f_C$ with content not exceeding $\mu$, one can find a
matrix $A \in \Theta(n,r)$ which defines that monomial. This proves
that the first basis given in \cite[(2.7)]{DG:PSA} coincides with the
basis described in \eqref{PSA:b} above.  One can argue similarly that
the second basis in \cite[(2.7)]{DG:PSA} coincides with the basis
\eqref{PSA:a}.\footnote{There is an error in the definition of the
second basis in \cite[(2.7)]{DG:PSA} (and similarly in
\cite[(3.9)]{DG:PSA}): one needs a notion of content {\em opposite} to
the one used to define the first basis there.}

By commutation relations given in \cite[Proposition 4.5]{DG:PSA} one
can rewrite a given basis element of the form \eqref{PSA:a} in the
form
\begin{equation} \label{PSA:c}
1_{\row(A)}\, ( \prod_{i<j}\divided{f_{ij}}{a_{ji}} 
\prod_{i<j} \divided{e_{ij}}{a_{ij}} )\, 1_{\col(A)}
\quad (A \in \Theta(n,r));
\end{equation}
similarly one can rewrite a given basis element of the form
\eqref{PSA:b} in the form
\begin{equation}\label{PSA:d}
1_{\row(A)}\, ( \prod_{i<j} \divided{e_{ij}}{a_{ij}}
\prod_{i<j} \divided{f_{ij}}{a_{ji}} )\, 1_{\col(A)} 
\quad (A \in \Theta(n,r))
\end{equation}
where $\row(A)$ (resp., $\col(A)$) is the vector of row (column) sums
in the matrix $A$.

Identify $S(\lambda)$ with $1_\lambda S(n,r) 1_\lambda$ and write
$S_\Z(\lambda) = 1_\lambda S_\Z(n,r) 1_\lambda$ and
${}_\lambda(S_\Z)_\mu = 1_\lambda S_\Z(n,r)1_\mu$.  Then
${}_\lambda(S_K)_\mu \simeq K \otimes_\Z {}_\lambda(S_\Z)_\mu$. The
following result provides a basis for $S(\lambda)$ specialized to any
field $K$.

\begin{thm}\label{thm:Zbas}
For $\lambda,\mu \in \Lambda(n,r)$ the bimodule ${}_\lambda(S_\Z)_\mu$
has a $\Z$-basis consisting of all elements
\begin{equation}\label{Zbas:a}
1_{\lambda}\, (\prod_{i<j} \divided{f_{ij}}{a_{ji}} 
\prod_{i<j} \divided{e_{ij}}{a_{ij}})\, 1_{\mu} \quad 
(A \in{}_\lambda \Theta_\mu)
\end{equation}
and another such basis consisting of all
\begin{equation}\label{Zbas:b}
1_{\lambda}\, (\prod_{i<j} \divided{e_{ij}}{a_{ij}} 
\prod_{i<j} \divided{f_{ij}}{a_{ji}})\, 1_{\mu} \quad 
(A \in{}_\lambda \Theta_\mu)
\end{equation}
where ${}_\lambda \Theta_\mu$ is the set of all $n \times n$ matrices
over $\N$ with row and column sums equal to $\lambda$, $\mu$, resp.,
and the products of $f$'s (resp., $e$'s) is taken with respect to some
fixed, but arbitrary, order. 
\end{thm}

\begin{proof}
This follows immediately from \eqref{PSA:c}, \eqref{PSA:d} above.
\end{proof}

\part{Generic Hecke algebras}

This part is concerned mainly with certain infinite-dimensional
quotients of a subalgebra $\U[0]$ of the enveloping algebra $\U =
\U(\gl_n)$, closely related to a modified form $\dot{\U}$ of $\U$.
The construction of $\dot{\U}$ was given by Lusztig in the quantum
case; its definition here, in \S\ref{sec:mf}, is virtually the
same. In particular, $\dot{\U}$ is an algebra over $\Q$ (without 1)
which is by definition the direct sum of various quotient spaces
${}_\lambda \U_\mu$ of $\U$. There is an infinite family
$\{1_\lambda\}_{\lambda \in \Z^n}$ of pairwise orthogonal idempotents
in $\dot{\U}$ which serves as a replacement for the identity, and
${}_\lambda \U_\mu = 1_\lambda \dot{\U}1_\mu$ for each $\lambda,\mu$.
We call $\dot{\U}(\lambda) := {}_\lambda \U_\lambda$ {\em generic}
Hecke algebras.  Each ${}_\lambda \U_\mu$ is a bimodule for
$\dot{\U}(\lambda)$-$\dot{\U}(\mu)$.

In \S\ref{sec:gen} we show that ${}_\lambda S_\mu$ and $S(\lambda)$
are homomorphic images of ${}_\lambda \U_\mu$ and $\dot{\U}(\lambda) =
{}_\lambda \U_\lambda$, and we obtain explicit integral PBW
bases for ${}_\lambda \U_\mu$. As an application, we find (see
\ref{rmk:symm}) that the group algebra $\Z\Sigma_r$ is a homomorphic
image of $\U_\Z[0]$, a certain subalgebra of $\U_\Z$.  In
\S\ref{sec:modsl} we consider the modified form of $\U(\sl_n)$, and
show that the generic algebras and bimodules are the same as the ones
constructed in terms of $\dot{\U}(\gl_n)$.  It makes little difference
whether one works with $\dot{\U}(\gl_n)$ or $\dot{\U}(\sl_n)$. We show
(see \S\ref{sec:cell}) that $\dot{\U}_\Z(\lambda)$, $S_\Z(\lambda)$
are cellular algebras in the sense of Graham and Lehrer \cite{GL}.
Finally, we describe an indexing set for the simple
$\dot{\U}(\lambda)$-modules (see \S\ref{sec:UKostka}).

\section{Lusztig's modified form of $\U(\gl_n)$}\label{sec:mf}

We work over $\Q$ in this section. We consider the modified form of
$\U = \U(\gl_n)$.

\subsection{}\label{Udot} 
Recall that $\U$ is the associative algebra with $1$ given by
generators $e_i$, $f_i$ ($1\le i \le n-1$), and $H_i$ ($1\le i \le n$)
subject to the usual relations given by the Lie algebra structure 
(see e.g.\ \cite[relations (R1)--(R5)]{DG:PSA}).

For $\lambda, \mu \in \Z^n$ we define ${}_\lambda \U_\mu$ (as a
vector space) to be the following quotient:
\begin{equation}
{}_\lambda \U_\mu = \textstyle
\U/(\sum_i (H_i-\lambda_i)\U + \sum_i \U(H_i-\mu_i)).
\end{equation}
Let $\pi_{\lambda,\mu}: \U \to {}_\lambda \U_\mu$ be the
canonical projection, and set $\dot{\U} = \oplus_{\lambda,\mu}
({}_\lambda \U_\mu)$.

Set $\alpha_i = \varepsilon_i - \varepsilon_{i+1}$ for $1\le i \le
n-1$, where $\{\varepsilon_1,\dots, \varepsilon_n\}$ is the canonical
basis of $\Z^n$. Consider the grading on $\U$ defined by putting $H_i \in
\U[0]$, $e_i \in \U[\alpha_i]$, $f_i \in \U[-\alpha_i]$ subject to
the requirement $\U[\nu']\U[\nu''] \subset \U[\nu'+\nu'']$. Then
$\U=\oplus_\nu \U[\nu]$, where $\nu$ runs over the set $\sum
\Z\alpha_i$ (the root lattice). We note that ${}_\lambda\U[\nu]_\mu :=
\pi_{\lambda,\mu} \U[\nu]$ is zero unless $\lambda-\mu=\nu$.

The above grading allows one to define a natural associative
$\Q$-algebra structure on $\dot{\U}$ inherited from that of $\U$. For
any $\lambda,\mu,\lambda',\mu' \in \Z^n$ and any $t\in
\U[\lambda-\mu]$, $s\in \U[\lambda'-\mu']$, the product
$\pi_{\lambda,\mu}(t) \pi_{\lambda',\mu'}(s)$ is defined to be equal
to $\pi_{\lambda,\mu'}(ts)$ if $\mu=\lambda'$ and is zero otherwise.

The elements $1_\lambda = \pi_{\lambda,\lambda}(1)$ satisfy the
relation $1_\lambda 1_\mu = \delta_{\lambda,\mu} 1_\lambda$, i.e.,
$\{ 1_\lambda\}$ is a family of orthogonal idempotents in
$\dot{\U}$. Moreover,
\begin{equation}\label{Udot:a0}
{}_\lambda \U_\mu = 1_\lambda \dot{\U} 1_\mu \quad (\lambda,\mu \in
\Z^n).
\end{equation}
The algebra $\dot{\U}$ does not have $1$, since the infinite sum
$\sum 1_\lambda$ does not belong to $\dot{\U}$.

\subsection{}\label{Udot:bi}
There is a natural $\U$-bimodule structure on $\dot{\U}$, defined by
the requirement
\begin{equation}\label{bi:a}
t \pi_{\lambda,\mu}(s) t' = \pi_{\lambda+\nu,\mu-\nu'}(tst')
\end{equation}
for all $t\in \U[\nu]$, $t'\in \U[\nu']$, $\lambda,\mu \in \Z^n$. 
Moreover, the following identities hold in the algebra $\dot{\U}$:
\begin{gather}\label{bi:b}
e_i 1_\lambda = 1_{\lambda+\alpha_i} e_i; \quad 
f_i 1_\lambda = 1_{\lambda-\alpha_i} f_i;\\
(e_if_j-f_je_i)1_\lambda = \delta_{ij} \lambda_i 1_\lambda 
\end{gather}
for all $\lambda \in \Z^n$, all $i,j$.  One can show that the algebra
$\dot{\U}$ is generated by all elements of the form
\begin{equation}\label{Udot:a}
e_i 1_\lambda,\ \ f_i 1_\lambda \quad (1 \le i \le n-1,\ \lambda \in
\Z^n).
\end{equation}

The algebra $\dot{\U}$ inherits a ``comultiplication'' from the
comultiplication on $\U$ (follow \cite[23.1.5]{Lusztig:book}) but we
shall not need it here.

\subsection{} \label{Udot:mod}
Following Lusztig \cite[23.1.4]{Lusztig:book} we say that a
$\dot{\U}$-module $M$ is {\em unital} if: (a) for all $m\in M$ one
has $1_\lambda m = 0$ for all but finitely many $\lambda \in \Z^n$;
(b) for any $m\in M$ one has $\sum_{\lambda} 1_\lambda m = m$ (sum
over all $\lambda \in \Z^n$).

If $M$ is a unital $\dot{\U}$-module, then one may regard it as a
$\U$-module with weight space decomposition, as follows. The weight
space decomposition $M = \oplus\, ({}_\lambda M)$ is given by setting
${}_\lambda M = 1_\lambda M$ and the action of $u\in \U$ is defined by
$u\cdot m = (u1_\lambda)m$ for any $m\in {}_\lambda M$, where
$u1_\lambda$ is regarded as an element of $\dot{\U}$ as in
\ref{Udot:bi}.

From the above one sees that the category of unital $\dot{\U}$-modules
is the same as the category of $\U$-modules which admit weight space
decompositions.

\subsection{} \label{UdotZ}
Let $\dot{\U}_\Z$ be the subalgebra of $\dot{\U}$ generated by all
\begin{equation}\label{UdotZ:a}
\divided{e_i}{c} 1_\lambda,\ \ \divided{f_i}{a} 1_\lambda  \quad 
(1 \le i \le n-1,\ \lambda \in \Z^n,\ a,c \ge 0).
\end{equation}
Let $B^+$, $B^-$ be any $\Z$-bases of $\U^+_\Z$, $\U^-_\Z$, resp. By
the analogue of \cite[23.2.1]{Lusztig:book} we have a $\Z$-basis for
$\dot{\U}_\Z$ consisting of all elements of the form
\begin{equation}\label{tri:1}
b^+ 1_\lambda b'^- \qquad (b^+ \in B^+,\ b'^- \in B^-,\ \lambda \in \Z^n)
\end{equation}
and another such basis consisting of the elements of the form 
\begin{equation}\label{tri:2}
b'^- 1_\lambda b^+ \qquad (b^+ \in B^+,\ b'^- \in B^-,\ \lambda \in \Z^n).
\end{equation}
It follows that there is a $\Z$-basis for $\dot{\U}_\Z$ consisting of
all elements of the form
\begin{equation}\label{UdotZ:c}
(\prod_{i<j} \divided{f_{ij}}{a_{ji}})\, 1_{\lambda^-(A)}\,
(\prod_{i<j} \divided{e_{ij}}{a_{ij}}) \quad (A \in \widetilde{\Theta}(n))
\end{equation}
and another $\Z$-basis consisting of all elements of the form
\begin{equation}\label{UdotZ:d}
(\prod_{i<j} \divided{e_{ij}}{a_{ij}})\, 1_{\lambda^+(A)}\,
(\prod_{i<j} \divided{f_{ij}}{a_{ji}}) \quad (A \in \widetilde{\Theta}(n))
\end{equation}
where $\widetilde{\Theta}(n)$ is the set of $n \times n$ integral
matrices with off-diagonal entries $\ge 0$, and where as usual
products of $f$'s (resp., $e$'s) are taken with respect to some fixed
order. The definition of $\lambda^+(A)$, $\lambda^-(A)$ here is the
same as in \ref{PSA}.

\section{The generic algebra $\dot{\U}(\lambda)={}_\lambda \U_\lambda$}
\label{sec:gen}

\subsection{}
From the definition of the multiplication in $\dot{\U}$ it is clear
that ${}_\lambda \U_\lambda = 1_\lambda \dot{\U} 1_\lambda$ is an
algebra and ${}_\lambda \U_\mu = 1_\lambda \dot{\U} 1_\mu$ is a
bimodule (for ${}_\lambda \U_\lambda$-${}_\mu \U_\mu$).  By
definition, $\dot{\U}$ is the direct sum of these bimodules, as
$\lambda,\mu$ range over $\Z^n$. We call the algebra
$\dot{\U}(\lambda) = {}_\lambda \U_\lambda$ a generic algebra.

We have the following analogue of Lemma \ref{lem:iso}, which shows
that one needs to consider $\dot{\U}(\lambda)$ and ${}_\lambda\U_\mu$
only in case $\lambda,\mu \in \Z^n$ are dominant. 

\begin{lem}\label{lem:Uiso} 
Let $\lambda, \mu \in \Z^n$.

\noindent(a) For any $w \in \Sigma_n$,
one has an algebra isomorphism $\dot{\U}(\lambda) \simeq
\dot{\U}(w \lambda)$.

\noindent(b) For any $w,w' \in \Sigma_n$, one has an isomorphism
${}_\lambda\U_\mu \simeq {}_{w\lambda}\U_{w' \mu}$ of
$\dot{\U}(\lambda)$-$\dot{\U}(\mu)$ bimodules.
\end{lem}

\begin{proof}
The Weyl group $W \simeq \Sigma_n$ acts as permutations on the Lie
algebra $\gl_n$ via the rule: $w e_{ij} = e_{w^{-1}(i), w^{-1}(j)}$.
(Here the $e_{ij}$ are the matrix units defined in \ref{Ubas}.) This
induces a corresponding action of $\Sigma_n$ on $\U$. Using this
action, the statements are easily checked. 
\end{proof}

\begin{prop}\label{prop:trunc}
For any $\lambda,\mu \in \Lambda(n,r)$ the linear map
$\psi_{\lambda,\mu}: \U \to 1_\lambda S(n,r) 1_\mu$ given by $u \to
1_\lambda \pi_{n,r}(u) 1_\mu$ is surjective, where $\pi_{n,r}: \U \to
S(n,r)$ is the map defined in \eqref{Ubas:d}.  Moreover, the
above linear map induces a linear surjection ${}_\lambda \U_\mu \to
1_\lambda S(n,r)1_\mu$. 
\end{prop}

\begin{proof}
By \cite[Proposition 4.3(a)]{DG:PSA}, in $S(n,r)$ we have the relation
$H_i 1_\mu = 1_\mu H_i = \mu_i1_\mu$, for any $i$ and any $\mu \in
\Z^n$. It follows that $H_i-\lambda_i$ acts on the left on $1_\lambda
S(n,r)1_\mu$ as zero; similarly, $H_i-\mu_i$ acts as zero on the
right.  Thus the map $\psi_{\lambda,\mu}$ sends any element of $\sum_i
(H_i-\lambda_i)\U + \sum_i \U(H_i-\mu_i)$ to zero, so its kernel
contains the kernel of ${}_\lambda \U_\mu$, and so
$\psi_{\lambda,\mu}$ factors through ${}_\lambda \U_\mu$. Moreover,
$\psi_{\lambda,\mu}$ is surjective since the map
$$
u \to \sum_{\lambda,\mu \in \Lambda(n,r)} 1_\lambda \pi_{n,r}(u) 1_\mu 
$$
is the surjection $\U \to S(n,r)$ of \eqref{Ubas:d}. All the
assertions of the proposition are now clear. 
\end{proof}

\begin{cor}
(a) The restriction to $\U_\Z$ of $\psi_{\lambda,\mu}$ 
surjects onto $1_\lambda S_\Z(n,r)1_\mu$. 

\noindent(b) For $\lambda\in \Lambda(n,r)$, $S(\lambda)$ is a
homomorphic image of $\dot{\U}(\lambda)$.  Similarly, $S_\Z(\lambda)$
is a homomorphic image of $\dot{\U}_\Z(\lambda)$.
\end{cor}

\begin{rmks}\label{rmk:symm}
(a) The map $\psi_{\lambda,\lambda}: \U \to S(\lambda)$ is not in
general an algebra map. (For instance, the product $f_ie_i$ will map
to something nonzero but $e_i$, $f_i$ themselves map to zero.)
However, the restriction of $\psi_{\lambda,\lambda}$ to the subalgebra
$\U[0]$ (see \ref{Udot}) of $\U$ will be an algebra map. So, for
$\lambda \in \Lambda(n,r)$, $S(\lambda)$ is a homomorphic image of
$\U[0]$ and $S_\Z(\lambda)$ is a homomorphic image of $\U_\Z[0]$.

(b) In particular, this shows that $\Z\Sigma_r$ is a quotient of
$\U_\Z[0]$; in other words, $\Z\Sigma_r$ is a subquotient of $\U_\Z$.
Similarly, $\Z\Sigma_r$ is a subquotient of $\dot{\U}_\Z$, since it is
a homomorphic image of ${}_\lambda(\U_\Z)_\lambda =
\dot{\U}_\Z(\lambda)$.  Thus for any field $K$ the group algebra
$K\Sigma_r$ is a subquotient of the hyperalgebra $\U_K=K\otimes_\Z
\U_\Z$ (and of $\dot{\U}_K=K\otimes_\Z \dot{\U}_\Z$).
\end{rmks}

The following result provides a basis for $\dot{\U}_K(\lambda):=
\dot{\U}_\Z(\lambda) \otimes_\Z K$ under specialization to any field
$K$.

\begin{prop}\label{thm:dotUbas}
For any $\lambda,\mu \in \Z^n$, a $\Z$-basis for the bimodule
${}_\lambda(\dot{\U}_\Z)_\mu = 1_\lambda \dot{\U}_\Z 1_\mu$ is given
by the set of all elements of the form
\begin{equation}\label{dotUbas:a}
1_{\lambda}\, (\prod_{i<j} \divided{f_{ij}}{a_{ji}} 
\prod_{i<j} \divided{e_{ij}}{a_{ij}})\, 1_{\mu} \quad 
(A \in {}_\lambda \widetilde{\Theta}_\mu)
\end{equation}
and another such basis is given by all elements of the form
\begin{equation}\label{dotUbas:b}
1_{\lambda}\, (\prod_{i<j} \divided{e_{ij}}{a_{ij}} 
\prod_{i<j} \divided{f_{ij}}{a_{ji}})\, 1_{\mu} \quad 
(A \in {}_\lambda \widetilde{\Theta}_\mu)
\end{equation}
where ${}_\lambda \widetilde{\Theta}_\mu$ is the set of all matrices
$A\in \widetilde{\Theta}(n)$ with row and column sums equal to
$\lambda$, $\mu$, resp. As usual, products of $f$'s (resp., $e$'s) is
taken with respect to some fixed, but arbitrary, order. 
\end{prop}

\begin{proof} 
This follows from \eqref{UdotZ:c}, \eqref{UdotZ:d} and commutation
formulas \eqref{bi:b}. The argument is similar to the argument for
\ref{thm:Zbas}.
\end{proof}

\begin{example}
We consider the case $n=2$. Let $\lambda = (\lambda_1,\lambda_2) \in
\Z^2$. Then the basis of $\dot{\U}_\Z(\lambda)$ described in
\eqref{dotUbas:a} consists of all elements of the form
\begin{equation}\label{example:a}
1_\lambda \divided{f}{a} \divided{e}{a} 1_\lambda \quad(a \ge 0)
\end{equation}
and the basis described in \eqref{dotUbas:b} consists of all elements of
the form
\begin{equation}\label{example:b}
1_\lambda \divided{e}{a} \divided{f}{a} 1_\lambda \quad(a \ge 0).
\end{equation}
If $\lambda \in \Lambda(2,r)$, the nonzero images in $S_\Z(\lambda)$
of the elements \eqref{example:a} give the following set of elements
\begin{equation}\label{example:c}
1_\lambda \divided{f}{a} \divided{e}{a} 1_\lambda \quad(0 \le a \le
\min(\lambda_1,\lambda_2))
\end{equation}
which was described in \eqref{Zbas:a}; the nonzero images in
$S_\Z(\lambda)$ of the elements \eqref{example:b} give the following
set of elements
\begin{equation}\label{example:d}
1_\lambda \divided{e}{a} \divided{f}{a} 1_\lambda \quad(0 \le a \le
\min(\lambda_1,\lambda_2))
\end{equation}
which was described in \eqref{Zbas:b}. The sets in \eqref{example:c},
\eqref{example:d} are bases of $S_\Z(\lambda)$. Thus $\dim S(\lambda)
= 1+\min(\lambda_1,\lambda_2)$.

Note that from this description it follows that $\dot{\U}(\lambda)$
(for $\lambda \in \Z^2$) is generated by the element $1_\lambda fe
1_\lambda$. It is also generated by the element $1_\lambda ef
1_\lambda$. It follows that the algebra $\dot{\U}(\lambda)$ is
commutative.  The same statements apply to $S(\lambda)$ (for
$\lambda\in \Lambda(2,r)$).  This does not hold for partitions with
more than 2 parts, in general.
\end{example}

\section{The modified form of $\U(\sl_n)$}\label{sec:modsl}

\subsection{} 
Note that $\U(\sl_n)$ is the subalgebra of $\U=\U(\gl_n)$ generated by
all $e_i$, $f_i$, $h_i:= H_i-H_{i+1}$ ($1 \le i \le n-1$); see e.g.
\cite[p.\ 1909]{DG:PSA}, sentence following Theorem 2.1.  The
definition of $\dot{\U}(\sl_n)$ is nearly the same as the definition
of $\dot{\U}(\gl_n)$. For $\lambda, \mu \in \Z^{n-1}$ one defines
${}_\lambda \U(\sl_n)_\mu$ (as a vector space) to be the following
quotient space:
\begin{equation}
\textstyle \U(\sl_n) /(\sum_i
(h_i-\bil{\lambda}{\alpha_i^\vee})\U(\sl_n) + \sum_i
\U(\sl_n)(h_i-\bil{\mu}{\alpha_i^\vee})).
\end{equation}
Set $\dot{\U}(\sl_n) = \oplus_{\lambda,\mu} ({}_\lambda
\U(\sl_n)_\mu)$. The rest of the construction is exactly the same as
for $\dot{\U}(\gl_n)$. The only difference between $\dot{\U}$ and
$\dot{\U}(\sl_n)$ is that the former has more idempotents than the
latter.

\subsection{} Given a weight $\lambda \in \Z^n$ (for $\gl_n$) we obtain
a corresponding weight $\widetilde{\lambda} \in \Z^{n-1}$ (for
$\sl_n$) as follows:
$$
\lambda \to \widetilde{\lambda}:= (\lambda_1-\lambda_2, 
\lambda_2-\lambda_3, \dots, \lambda_{n-1}-\lambda_n).
$$

The restriction to $\U(\sl_n)$ of the map $\pi_{n,r}$ is still a
surjection onto $S_\Q(n,r)$. This follows from the decomposition
$\gl_n = \sl_n \oplus \Q I$; the image of $I$ in $\U(\gl_n)$ acts as
scalars.  Hence the Schur algebra $S_\Q(n,r)$ is a homomorphic image
of $\U(\sl_n)$.  This restriction is compatible with integral forms;
i.e., the restriction of $\pi_{n,r}$ to $\U_\Z(\sl_n)$ surjects onto
$S_\Z(n,r)$.

\begin{prop}
For $\lambda,\mu \in \Z^n$, the natural map ${}_{\widetilde{\lambda}}
\U_\Z(\sl_n)_{\widetilde{\mu}} \to {}_\lambda \U_\Z(\gl_n)_\mu$ is an
isomorphism of bimodules. In case $\mu=\lambda$ it is an isomorphism
of algebras. 
\end{prop}

\begin{proof}
One can adapt the argument for \cite[23.2.5]{Lusztig:book} to prove
the map is a $\Z$-linear isomorphism.  In fact, one can easily check
that there are $\Z$-bases of ${}_{\widetilde{\lambda}}
\U_\Z(\sl_n)_{\widetilde{\mu}}$ of the form \eqref{dotUbas:a} and
\eqref{dotUbas:b}, with $1_\lambda$, $1_\mu$ there replaced by
$1_{\widetilde{\lambda}}$, $1_{\widetilde{\mu}}$. The bijection is now
clear.  The last claim follows from the compatibility of the
isomorphism with the product in $\dot{\U}(\gl_n)$, $\dot{\U}(\sl_n)$.
\end{proof}

Thus we see that, for the study of the Hecke algebras $S(\lambda)$,
one can use a descent from generic subalgebras of either algebra
$\dot{\U}(\sl_n)$ or $\dot{\U}(\gl_n)$. The difference between these
viewpoints is merely notational.

\section{Cellularity of $\dot{\U}(\lambda)$, $S(\lambda)$} 
\label{sec:cell}

Lusztig \cite{Lusztig:CB} showed that the positive part $\UU^+$ of a
quantized enveloping algebra $\UU$ has a canonical basis.  In
\cite[Part IV]{Lusztig:book} the canonical basis is extended to
$\dot{\UU}$.  We show that the canonical basis on $\dot{\U}$ induces
compatible canonical bases on $\dot{\U}(\lambda)$, $S(\lambda)$.  We
apply this information to show that $\dot{\U}(\lambda)$ and
$S(\lambda)$ inherit canonical bases from the canonical basis on
$\dot{\U}$, and that these bases are cellular bases.  In particular,
this shows that $\dot{\U}(\lambda)$ and $S(\lambda)$ are cellular
algebras. (One can see the cellularity in other ways; see e.g.\
\cite[Theorem 6.6]{DJM}, \cite[Chapter 4, Exer.~13]{Mathas:book}.)

\subsection{} We recall from \cite{GL} the definition 
of cellularity. Let $A$ be an associative algebra over a ring $R$
(commutative with $1$). We do not insist that $A$ has $1$, nor do we
insist that $A$ be finite-dimensional. A {\em cell datum} for $A$ is a
quadruple $(\Lambda,M,C,\iota)$ where:

(a) $\Lambda$ (the set of weights) is partially ordered by $\ge$, $M$
is a function from $\Lambda$ to the class of finite sets, and $C$ is
an injective function
$$
 C: \coprod_{\lambda \in \Lambda} (M(\lambda) \times M(\lambda)) \to A
$$
with image an $R$-basis of $A$. If $\lambda \in \Lambda$ and $S,T \in
M(\lambda)$ then one writes $C^\lambda_{S,T}$ for $C(S,T)$. 

(b) The map $\iota$ is an $R$-linear involutory antiautomorphism of
$A$ such that $\iota(C^\lambda_{S,T}) = C^\lambda_{T,S}$.

(c) If $\lambda \in \Lambda$ and $S,T \in M(\lambda)$ then for all $a
\in A$ one has 
$$
a C^\lambda_{S,T} \equiv \sum_{S' \in M(\lambda)} r_a(S',S)
C^\lambda_{S',T}
\quad \mod{A( > \lambda)},
$$
where $r_a(S',S)$ is independent of $T$ and $A(> \lambda)$ is the
$R$-submodule of $A$ generated by the set of all $C^\mu_{S'',T''}$
such that $\mu > \lambda$, $S'', T'' \in M(\mu)$. 

\medskip

Any algebra $A$ that possesses a cell datum is said to be cellular and
the basis $\{ C^\lambda_{S,T} \}$ produced by its cell datum is its
cellular basis.

\begin{rmk}
For our applications, it is convenient to use the order on $\Lambda$
opposite to the order in the usual definition, so we have reversed
the usual ordering of weights in the definition above.
\end{rmk}

\subsection{} Following R.M.\ Green \cite{RG} we say that a cell 
datum $(\Lambda,M,C,\iota)$ is of {\em profinite} type if $\Lambda$ is
infinite and if for each $\lambda \in \Lambda$, the set
$
\{ \mu \in \Lambda: \mu \le \lambda  \} 
$
is finite. 

\medskip

The following useful lemma follows immediately from the definitions. 

\begin{lem}\label{lem:cell}
(a) Let $(\Lambda, M, C, \iota)$ be a cell datum for $A$.  Let $e \in
A$ be an idempotent fixed by the involution $\iota$. Then $eAe$ is
cellular, with cell datum $(\overline{\Lambda}, \overline{M},
\overline{C}, \overline{\iota})$, where
$\overline{\Lambda} = \{ \lambda \in \Lambda: eC^\lambda_{S,T}e \ne
0, \text{some } S,T \in M(\lambda) \}$; 
$\overline{M}(\lambda) = \{ S,T \in M(\lambda) : e C^\lambda_{S,T} e 
\ne 0 \}$ for any $\lambda \in \overline{\Lambda}$; 
$\overline{C}$ is defined by 
$\overline{C}^\lambda_{S,T} = e C^\lambda_{S,T} e$ whenever $\lambda
\in \overline{\Lambda}$, $S,T \in \overline{M}(\lambda)$; and
$\overline{\iota}$ is the restriction of $\iota$ to $eAe$.

\noindent(b) If the original cell datum $(\Lambda, M, C, \iota)$ is of
profinite type, then $(\overline{\Lambda}, \overline{M}, \overline{C},
\overline{\iota})$ is also of profinite type, provided
$\overline{\Lambda}$ is infinite.
\end{lem}

\subsection{}\label{base}
The remarks in this subsection apply to a quantized enveloping algebra
$\UU = \UU(\g)$ correponding to an arbitrary reductive Lie algebra
$\g$.  We follow the setup and notation of \cite{Lusztig:book}.  Let
$\U = \U(\g)$.  Lusztig \cite{Lusztig:book} shows that the canonical
basis on $\UU^+$ can be extended to a canonical basis of the modified
form $\dot{\UU}$.  Moreover, in \cite[Chapter 29]{Lusztig:book}
Lusztig shows that the canonical basis of $\dot{\UU}$ is a cellular
basis (of profinite type). This is spelled out in greater detail in
\cite[2.5]{D:gen}.

In \cite{Lusztig:CB} it was pointed out that the canonical basis on
$\UU^+$ (which is an $\A$-basis for $\UU^+_\A$) corresponds under
specialization to a $\Z$-basis of the plus part $\UU^+_1 = \U^+_\Z$ of
$\U_\Z$.  By \cite[23.2]{Lusztig:book}, one has an
$\A$-basis of $\dot{\UU}$ consisting of elements of the form
\begin{equation}\label{base:1}
b^+ 1_\lambda b^{\prime-}.
\end{equation}
Similarly one has another $\A$-basis consisting of elements of the
form
\begin{equation}\label{base:2}
b^{\prime-}1_\lambda b^+ .
\end{equation}
In both sets of elements above, $b^+$, (resp., $b^{\prime-}$) vary
independently over any $\A$-basis of $\UU^+_\A$, (resp., $\UU^-_\A$).
It follows that specializing $v$ to $1$ takes $\dot{\UU}_\A$ to
$\dot{\U}_\Z$.  Moreover, any $\A$-basis for $\dot{\UU}_\A$ of the
form \eqref{base:1} or \eqref{base:2} will correspond to a $\Z$-basis
of $\dot{\U}_\Z$. In particular, the canonical basis on $\dot{\UU}_\A$
corresponds under specialization to a $\Z$-basis of $\dot{\U}_\Z$; we
call this the canonical basis of $\dot{\U}$. It is a cellular basis of
$\dot{\U}_\Z$ (of profinite type).

\medskip

In the rest of this section, $\U = \U(\gl_n)$ or $\U(\sl_n)$, and
we set $X = \Z^n$ or $\Z^{n-1}$. 

\begin{lem}\label{lem:CB4S}
The image of the canonical basis under the quotient map 
$\dot{\U}_\Z \to S_\Z(n,r)$ is a $\Z$-basis of $S_\Z(n,r)$.  
(We call it the canonical basis of $S_\Z(n,r)$).
\end{lem}

\begin{proof}
By Donkin \cite{Donkin:SA1,Donkin:SA2}, $S_\Z(n,r)$ is a generalized
Schur algebra, defined by the saturated set $\pi =
\Pi^+(\E_\Z^{\otimes r})$ (the set of dominant weights occuring in the
tensor space $\E_\Z^{\otimes r}$).  Since $S_\Z(n,r) \simeq
\dot{\U}/\dot{\U}[P]$, where $P$ is the complement of $\pi$ in the set
of dominant weights, the claim follows by the analogue of
\cite[29.2]{Lusztig:book}.
\end{proof}

\begin{rmk}
We have chosen, for simplicity, to obtain the canonical basis on
$S(n,r)$ by descent from the canonical basis of $\dot{\U}$. In
\cite{Du} this is approached (for the $q$-Schur algebra) the other way
around, by building up from the Kahzdan-Lusztig basis for the Hecke
algebra of type $A$.  See \cite{BLM} for yet another approach.
\end{rmk}

\begin{thm} \label{thm:cellular}
(a) For $\lambda \in X$, the generic algebra $\dot{\U}_\Z(\lambda)$ is
a cellular subalgebra of $\dot{\U}$, with cell datum of procellular
type. The cellular basis determined by the cell datum is inherited
from the canonical basis of $\dot{\U}_\Z$; we call this cellular basis
the canonical basis of $\dot{\U}_\Z(\lambda)$. 
 
\noindent(b) For $\lambda \in \Lambda(n,r)$, the Hecke algebra
$S_\Z(\lambda)$ is a cellular subalgebra of the Schur algebra
$S_\Z(n,r)$. The cellular basis determined by its cell datum is
inherited from the canonical basis of $S_\Z(n,r)$; we call this 
cellular basis the canonical basis of $S_\Z(\lambda)$. 

\noindent(c) For $\lambda \in \Lambda(n,r)$, the quotient map $\U \to
S(n,r)$ maps the canonical basis of $\dot{\U}_\Z(\lambda)$ (or of
$\dot{\U}_\Z(\widetilde{\lambda})$) onto the canonical basis of
$S_\Z(\lambda)$.
\end{thm}

\begin{proof}
This follows immediately from Lemmas \ref{lem:cell} and
\ref{lem:CB4S}, since the idempotent $1_\lambda$ is fixed by the
involution.
\end{proof}

\section{Simple $\dot{\U}_K(\lambda)$-modules} \label{sec:UKostka}

We work over an arbitrary infinite field $K$ in this section.  Set
$\dot{\U}_K = K\otimes_\Z \dot{\U}_\Z$, and similarly for
$\dot{\U}_K(\lambda)$, ${}_\lambda(\U_K)_\mu$. The simple unital
$\dot{\U}_K$-modules are the same as the simple rational
$\GL_n(K)$-modules. This follows, for instance, from the $v=1$
analogue of \cite[Chapter 31]{Lusztig:book}.
(Alternatively, see the discussion in \cite[II, 1.20]{Jantzen}.)

\subsection{}
Let $\lambda \in \Z^n$ be given. Let $\calF_\lambda$ be the functor $V
\to 1_\lambda V$ from unital $\dot{\U}_K$-modules to left
$\dot{\U}_K(\lambda)$-modules. This is an exact covariant functor
mapping simple modules to simple modules or $0$.  Let $L(\mu)$ be the
simple $\dot{\U}_K$-module corresponding to a dominant weight $\mu \in
\Z^n$.  By \cite[(6.2g)]{Green:book}, the collection of all nonzero
$\calF_\lambda L(\mu)$ forms a complete set of simple
$\dot{\U}_K(\lambda)$-modules.

\subsection{}
We write $\GL_n = \GL_n(K)$.  For a dominant $\mu \in \Z^n$, let
$I(\mu)$ denote the injective envelope of $L(\mu)$ in the category of
rational $\GL_n$-modules. Since induction takes injectives to
injectives, the module $\ind_{\T_n}^{\GL_n} K_\lambda$ is injective as
rational $\GL_n$-module, for any $\lambda \in \Z^n$.  Write $(
\ind_{\T_n}^{\GL_n} K_\lambda : I(\mu) )$ for the multiplicity of
$I(\mu)$ in a Krull-Schmidt decomposition of $\ind_{\T_n}^{\GL_n}
K_\lambda$. By Frobenius reciprocity one has an equality
\begin{equation}
(\ind_{\T_n}^{\GL_n} K_\lambda : I(\mu) ) = \dim_K 1_\lambda L(\mu)
\end{equation}
and this shows, in particular, that the multiplicities on the
left-hand-side are finite. Let us denote the integer in the equality 
by $\mathbf{K}^\prime_{\mu\lambda}$. The following analogue of
\ref{prop:Kd} is now clear.

\begin{prop} Let $\lambda \in \Z^n$ be fixed.
The isomorphism classes of simple $\dot{\U}_K(\lambda)$-modules are
the ${}_\lambda L(\mu) = 1_\lambda L(\mu)$ for which
$\mathbf{K}^\prime_{\mu\lambda} \ne 0$.
\end{prop}

\bibliographystyle{amsalpha}

\end{document}